\documentclass{article}

\usepackage[english]{babel}
\usepackage{textcomp}
\usepackage{amscd}
\usepackage{amssymb}
\usepackage{amsmath}
\usepackage{graphicx}
\usepackage{tikz}
\usepackage{ifthen}
\usetikzlibrary{snakes,arrows,shapes}

\newtheorem{proposition}{Proposition}

\newcommand{\Ar}{\mathrm{Area}}

\begin{document}
\DeclareGraphicsExtensions{.jpg,.pdf,.png}

\begin{center}
{\huge \textbf{Polymetric brick wall patterns}}

\smallskip

{\huge \textbf{and two-dimensional}}

\medskip

{\huge \textbf{ substitutions}}

\medskip

{ \large Michel Dekking}\\

\bigskip

3TU Applied Mathematics Institute and Delft
University of Technology, Faculty EWI, P.O.~Box 5031, 2600 GA
Delft, The Netherlands.

E-mail: f.m.dekking@tudelft.nl

\end{center}

\bigskip

\begin{abstract} Polymetric walls are walls built from bricks in more than one size. Architects and builders want to built polymetric walls that satisfy certain structural and aesthetical constraints. In a recent paper by de Jong,  Vindu$\check{s}$ka, Hans and Post these problems are solved by integer programming techniques, which can be very time consuming for patterns consisting of more than 40 bricks. Here we give an extremely fast method, generating patterns of arbitrary size.
\end{abstract}


\section{Introduction}

How to built a brick wall with different sizes of bricks? This problem was recently discussed in the paper \cite{JVHP}.

\begin{figure}[h!]%
\begin{center}
\includegraphics[height=4cm]{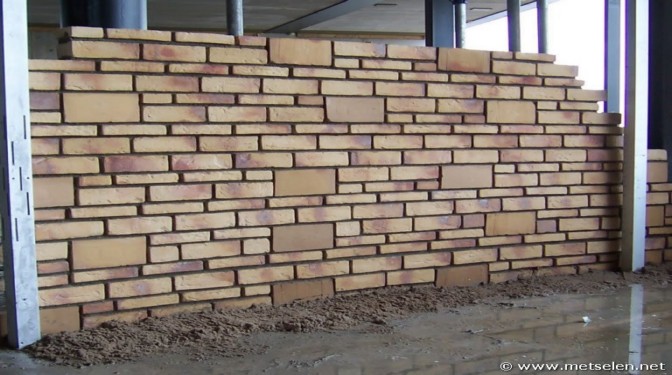}\\  %
  \caption{\footnotesize An example of a polymetric brick wall from www.metselen.net/fotoalbum.html}\label{Fig:1}
\end{center}
\end{figure}


Let us try to formulate  the requirements put forward by architects and builders. These are
\begin{itemize}
\item[[R1\!\!\!]] The absence of regularity in the patterns, as e.g., periodicity
\item[[R2\!\!\!]] Long horizontal joints, to facilitate the building process
\item[[R3\!\!\!]] The absence of long vertical joints, for sturdiness and strength
\item[[R4\!\!\!]] Pleasing aesthetics, as e.g., not too many large bricks close together
\end{itemize}

\noindent
One can appreciate the fulfilment of all four requirements in Figure~\ref{Fig:1}.

In \cite{JVHP} the four requirements are met by considering the construction of a brick wall pattern as a Pallet Loading Problem, which is then solved with a Mixed Integer Linear Program. Here we take a completely different approach, constructing the patterns as pseudo-self-similar tilings generated by 2-D substitutions. Well known properties of such tilings automatically yield the fulfilment of [R1], and also more or less [R4].
Requirement [R2] can very easily be built in, so that [R3] becomes the main challenge.

The style of this note will be somewhat informal, for precise definitions and more context see the papers \cite{Solom} and  \cite{Priebe-primer}. See Section 14 (and the book cover) of Allouche and Shallit's monograph \cite{Allou} for the subclass of tilings generated by multi-dimensional automatic sequences. Since then the literature on tilings has expanded  a lot (cf.~\cite{Priebe-primer}).

\section{Prouhet-Thue-Morse patterns}

 The Prouhet-Thue-Morse  sequence 0,1,1,0,1,0,0,1,... is the \emph{one}-dimensional
sequence  generated by the substitution $0\rightarrow 0\,1,\; 1\rightarrow 1\,0$. It is folklore that in the ``same way" one generates a 2D-sequence (see, e.g., \cite{Photonic}) by the 2D-substitution: \\
$$
\begin{tabular}{ccccccccc}
$ $ & $ $       & $1$ & $0$\; &$ $ & $ $ & $ $       & $0$ & $1$  \\
$0$ & $\mapsto$ & $0$ & $1$, &$ $ & $1$ & $\mapsto$ & $1$ & $0$.
\end{tabular}
$$
Usually this sequence is then represented in the plane by white and black squares. When we
 take  $2\times1$ and $3\times1$ bricks this corresponds to the following substitution rules

\begin{center}
\begin{tikzpicture}[scale=.2,rounded corners=0]
\def\bricktwoone(#1){\path (#1) coordinate (P0); \draw [fill=red!50!yellow, draw=gray, ultra thick]  (P0) rectangle +(2,1)}
\def\brickthreeone(#1){\path (#1) coordinate (P0);  \draw [fill=red!10!brown, draw=gray, ultra thick] (P0) rectangle +(3,1)}
\def\sigma(#1,#2){\path (#2) coordinate (P0);
\path (P0)++(0,0) coordinate (P210); \path (P0)++(2,0) coordinate  (P211);
\path (P0)++(0,1) coordinate (P212); \path (P0)++(3,1) coordinate  (P213);
\ifthenelse{#1=21}{\bricktwoone(P210);\brickthreeone(P211); \brickthreeone(P212);\bricktwoone(P213)}{};
\path (P0)++(0,0) coordinate (P310); \path (P0)++(3,0) coordinate (P311);
\path (P0)++(0,1) coordinate (P312); \path (P0)++(2,1) coordinate (P313); 
\ifthenelse{#1=31}{\brickthreeone(P310);\bricktwoone(P311);\bricktwoone(P312);\brickthreeone(P313)}{};
}; 
\path (0,0) coordinate (B21); \path (6,0) coordinate (SB21); \path (20,0) coordinate (B31); \path (26,0) coordinate (SB31);
\bricktwoone(B21); \sigma(21,SB21); \node at (4,0.4) {$\mapsto$};
\brickthreeone(B31); \sigma(31,SB31); \node at (24.7,0.4) {$\mapsto$};
\end{tikzpicture}
\end{center}

The third iteration of this substitution starting from the $2\times1$  brick $B_{21}$ yields the following pattern

\begin{figure}[h!]%
\begin{center}
 \includegraphics[height=1.5cm]{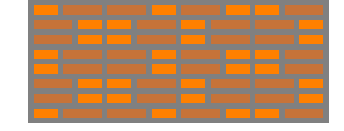}
 \end{center}
\end{figure}

Obviously requirement [R3] will  be violated, since we obtain vertical segments crossing the whole pattern.
To cope with this, we slightly change the previous substitution to the substitution $\tau$ given by

\begin{center}
\begin{tikzpicture}[scale=.2,rounded corners=0]
\def\bricktwoone(#1){\path (#1) coordinate (P0); \draw [fill=red!50!yellow, draw=gray, ultra thick]  (P0) rectangle +(2,1)}
\def\brickthreeone(#1){\path (#1) coordinate (P0);  \draw [fill=red!10!brown, draw=gray, ultra thick] (P0) rectangle +(3,1)}
\def\sigma(#1,#2){\path (#2) coordinate (P0);
\path (P0)++(0,0) coordinate (P210); \path (P0)++(2,0) coordinate  (P211);
\path (P0)++(1,1) coordinate (P212); \path (P0)++(4,1) coordinate  (P213);
\ifthenelse{#1=21}{\bricktwoone(P210);\brickthreeone(P211); \brickthreeone(P212);\bricktwoone(P213)}{};
\path (P0)++(0,0) coordinate (P310); \path (P0)++(3,0) coordinate (P311);
\path (P0)++(1,1) coordinate (P312); \path (P0)++(3,1) coordinate (P313); 
\ifthenelse{#1=31}{\brickthreeone(P310);\bricktwoone(P311);\bricktwoone(P312);\brickthreeone(P313)}{};
}; 
\path (0,0) coordinate (B21); \path (6,0) coordinate (SB21); \path (20,0) coordinate (B31); \path (26,0) coordinate (SB31);
\bricktwoone(B21); \sigma(21,SB21); \node at (4,0.4) {$\mapsto$};
\brickthreeone(B31); \sigma(31,SB31); \node at (24.7,0.4) {$\mapsto$};
\end{tikzpicture}
\end{center}

\begin{figure}[h!]%
\begin{center}
 \includegraphics[height=4cm]{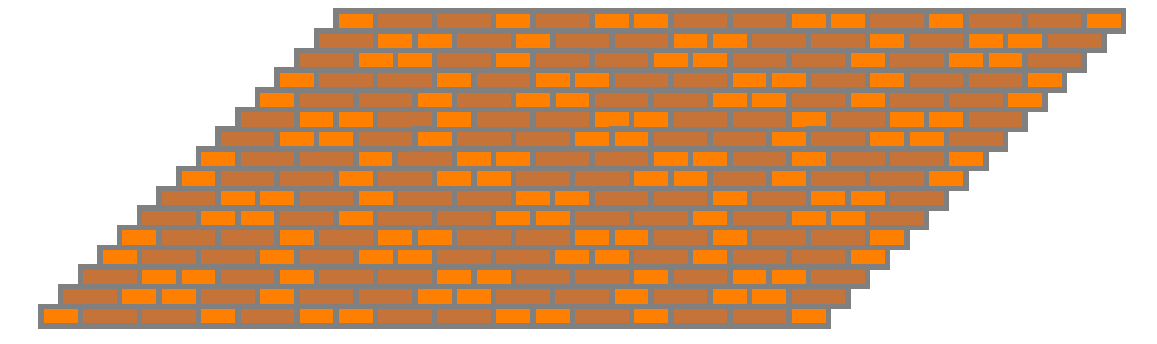}\\  %
  \caption{\footnotesize Skewed 2D-Prouhet-Thue-Morse: $\tau^4(B_{21})$.}\label{Fig:P-T-M-order-4}
\end{center}
\end{figure}

\noindent
Now a fourth order iteration of $\tau$ starting from the $2\times1$  brick yields the pattern in Figure~\ref{Fig:P-T-M-order-4}.

These skewed P-T-M patterns nicely satisfy [R3], see the following proposition.

\begin{proposition} For any $n$ the vertical segments in the patterns $\tau^n(B_{21})$ and $\tau^n(B_{31})$ have length at most 2.
\end{proposition}

\emph{Proof:} By induction on $n$, using  the following figure.
\begin{center}
\includegraphics[height=2cm]{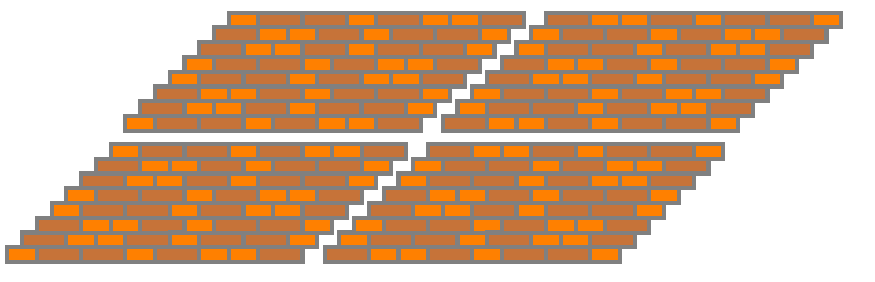}
\end{center}
\vspace*{-0.9cm}
\hfill $\square$

\section{General patterns}

We would like the pattern to be generated by a geometric 2D-substitution, mapping rectangles $i\times j$ (bricks $B_{ij}$)  to unions of rectangles (bricks), resulting in a self-similar tiling (see \cite{Priebe-primer} for the basic definitions). However, as already clear from the first P-T-M example above, such self-similar tilings will violate the bounded vertical segments condition [R3]. Thus we should  rather consider pseudo-self-similar or  pseudo-self-affine tilings as studied by Priebe and Solomyak (\cite{PriebeSolom}). In a self-affine (rectangular) tiling with expansion coefficients $\lambda_1$ and $\lambda_2$, all the bricks
(usually called proto-tiles in the literature) are scaled by a factor $\lambda_1$ in the $x$-direction and a factor $\lambda_2$ in the $y$-direction, and then are exactly tiled by translations of the original set of bricks. Note by the way that the P-T-M-substitution is self-similar with two $1\times 1$ bricks (expansions $\lambda_1=2$ and $\lambda_2=2$) , but not with a set consisting of the $2\times 1$ and a  $3\times 1$ brick.

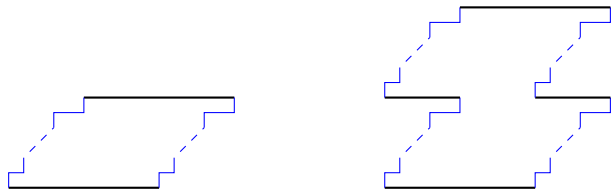
\begin{figure}[h!]%
\begin{center}
\begin{tikzpicture}[scale=.2,rounded corners=0]
\def\ypath(#1){\path (#1) coordinate (P0);
\path (P0)+(5,6) coordinate (Q0);
\path (P0)+(1.4,2.4) coordinate (D0);
\path (Q0)+(-2,-2) coordinate (D1);
\draw [blue] (P0)--++(0,1)--++(1,0)--++(0,1); \draw [blue] (Q0)--++(0,-1)--++(-2,0)--++(0,-1);
\draw [dashed, blue, thin] (D0)--(D1);}
\path (0,0) coordinate (left); \path (10,0) coordinate (right);
\path (5,6) coordinate (leftt); \path (15,6) coordinate (rightt);
\draw[thick] (left)--(right); \draw[thick] (leftt)--(rightt);
\ypath(left); \ypath(right);
\path (25,0) coordinate (shift);
\path (0,0)++(shift) coordinate (left); \path (10,0)++(shift) coordinate (right);
\path (5,6)++(shift) coordinate (leftt); \path (15,6)++(shift) coordinate (rightt);
\draw[thick] (left)--(right);
\ypath(left); \ypath(right);
\path (25,6) coordinate (shift);
\path (0,0)++(shift) coordinate (left); \path (10,0)++(shift) coordinate (right);
\path (5,6)++(shift) coordinate (lefttt); \path (15,6)++(shift) coordinate (righttt);
\draw[thick] (leftt)--(left); \draw[thick] (right)--(rightt);
\ypath(left); \ypath(right); \draw[thick] (lefttt)--(righttt);
\end{tikzpicture}
\end{center}
\caption{Left: Image region for a $i\times 1$ brick. Right: Image region for a $i\times 2$ brick.}\label{Fig:PP-prop}
\end{figure}

In a pseudo-self-affine tiling the bricks do not tile exactly the expanded bricks, but have bricks `sticking out' and `sticking in', in a way that still preserves the tiling property. We will consider a subclass of 2D-substitutions which have the desired properties, which we call of PP-type, i.e., which have the periodic parallelogram property. With this we mean that the image of a $i\times 1$ brick is a region as in Figure~\ref{Fig:PP-prop}, left. The black horizontal segments have length $i\lambda_1$, and the two curves on the left and the right are translated copies of each other. For a $i\times j$ brick this is repeated periodically $j$ times, cf. Figure~\ref{Fig:PP-prop}, right.

A simple example of a PP-type substitution with $\lambda_1=\lambda_2=2$ and the  three brick alphabet $\mathcal{A}=\{B_{11},B_{21},B_{22}\}$ is given by

\begin{center}
\begin{tikzpicture}[scale=.2,rounded corners=0]
\def\brickoneone(#1){\path (#1) coordinate (P0); \draw [fill=red!40!yellow, draw=gray, ultra thick]  (P0) rectangle +(1,1)}
\def\bricktwoone(#1){\path (#1) coordinate (P0); \draw [fill=red!60!yellow, draw=gray, ultra thick]  (P0) rectangle +(2,1)}
\def\bricktwotwo(#1){\path (#1) coordinate (P0);  \draw [fill=red!10!brown, draw=gray, ultra thick]  (P0) rectangle +(2,2)}
\def\sigma(#1,#2){\path (#2) coordinate (P0);
\path (P0)++(-1,0) coordinate (P110); \path (P0)++(0,1) coordinate  (P111);
\ifthenelse{#1=11}{\bricktwoone(P110);\bricktwoone(P111)}{};
\path (P0)++(-1,0) coordinate (P210); \path (P0)++(1,0) coordinate (P211);
\path (P0)++(0,1) coordinate (P212); \path (P0)++(3,1) coordinate (P213);
\ifthenelse{#1=21}{\bricktwoone(P210);\bricktwotwo(P211);\brickoneone(P212);\brickoneone(P213)}{};
\path (P0)++(-1,0) coordinate (P220); \path (P0)++(1,0) coordinate  (P221);
\path (P0)++(0,1) coordinate  (P222); \path (P0)++(3,1) coordinate  (P223);
\path (P0)++(-1,2) coordinate (P224); \path (P0)++(1,2) coordinate  (P225);
\path (P0)++(2,2) coordinate  (P226); \path (P0)++(0,3) coordinate  (P227);
\path (P0)++(2,3) coordinate  (P228);
\ifthenelse{#1=22}{\bricktwoone(P220);\bricktwotwo(P221);\brickoneone(P222);\brickoneone(P223);\bricktwoone(P224);\brickoneone(P225);\brickoneone(P226);%
\bricktwoone(P227);\bricktwoone(P228)}{};
}; 
\path (0,0) coordinate (B11); \path (6,0) coordinate (SB11);
\path (13,0) coordinate (B21); \path (20,0) coordinate (SB21);
\path (30,0) coordinate (B22); \path (37,0) coordinate (SB22);
\brickoneone(B11); \sigma(11,SB11); \node at (3,0.4) {$\mapsto$};
\bricktwoone(B21); \sigma(21,SB21); \node at (16.7,0.4) {$\mapsto$};
\bricktwotwo(B22); \sigma(22,SB22); \node at (33.7,0.8) {$\mapsto$};
\end{tikzpicture}
\end{center}
\noindent
The patterns $\sigma^3(B_{11}), \sigma^3(B_{21})$ and $\sigma^3(B_{22})$ are
\begin{center}
\includegraphics[height=2.5cm]{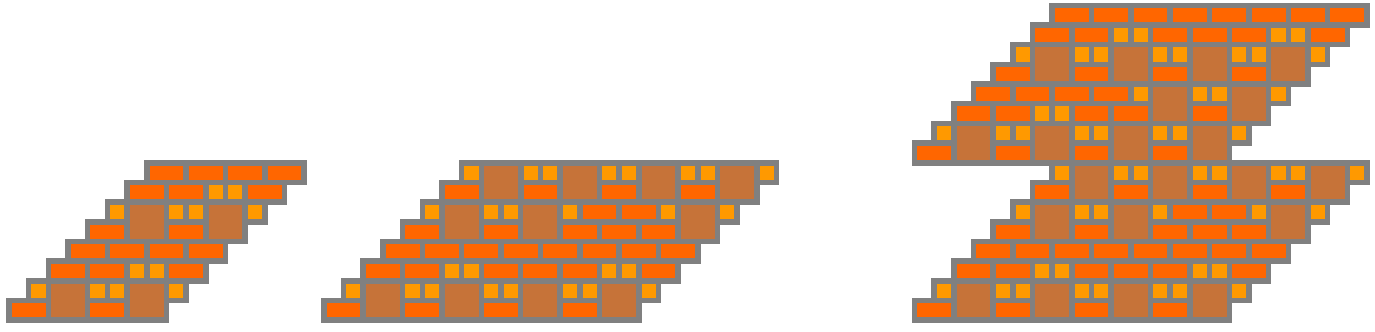}
\end{center}
The pattern $\sigma^4(B_{22})$ is given in Figure~\ref{Fig:max-11}.
\begin{figure}[h!]%
\begin{center}
 \includegraphics[height=6cm]{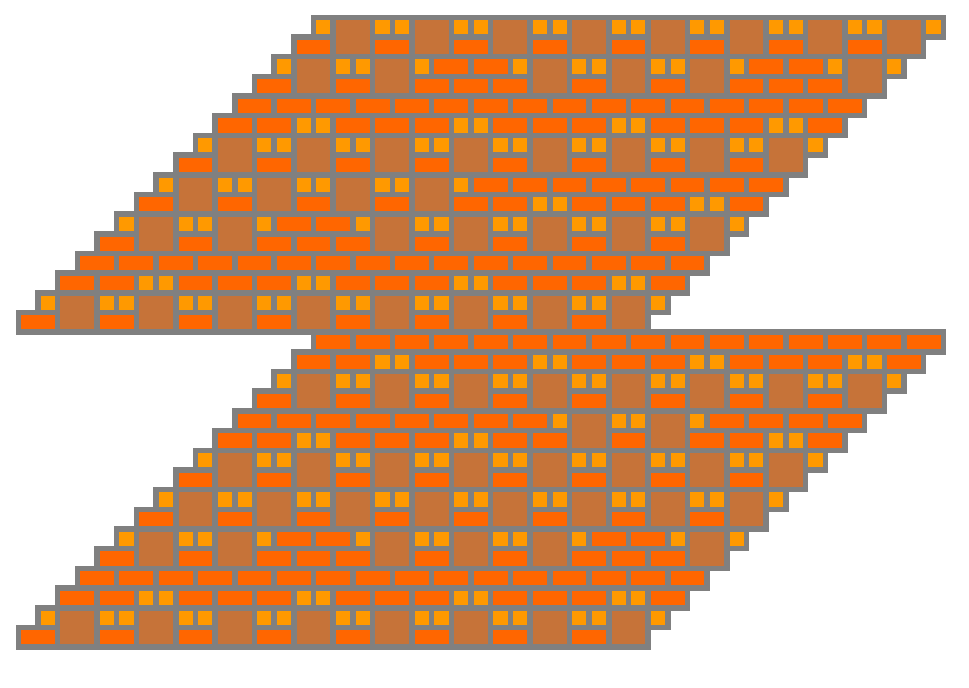}
  \caption{\footnotesize  $\sigma^4(B_{22})$.}\label{Fig:max-11}
\end{center}
\end{figure}

\medskip

Let $v_{\rm max}(ij,n)$ be the maximal length of a vertical segment in $\sigma^n(B_{ij})$ for $n\ge 1$ and $B_{ij}\in \mathcal{A}$.
For the substitution above $v_{\rm max}(11,1)=1, v_{\rm max}(21,1)=2, v_{\rm max}(22,1)=3$, and one can prove that for this substitution
$$v_{\rm max}:=\sup\{v_{\rm max}(ij,n):  B_{ij}\in \mathcal{A}, \;n\ge 1\}=v_{\rm max}(22,3)=11,$$
but in general it is not an easy task to determine the largest possible length $v_{\rm max}$ of a vertical segment  (which, of course, can be infinite).

The following approach generates patterns satisfying [R3] in an easy way. We say the pattern $\sigma(B_{ij})$ \emph{has no crossings} if there are no vertical segments consisting of boundaries of the bricks in $\sigma(B_{ij})$ that connect two points on the boundary of the region determined by $\sigma(B_{ij})$.
In the substitution above, for example, $\sigma(B_{11})$ and $\sigma(B_{22})$ do not have a crossing, but $\sigma(B_{21})$ does.

\begin{proposition}\label{Prop:nocross} Let $\sigma$ be a {\rm{PP}}-type substitution with expansions $\lambda_1$ and $\lambda_2$. Suppose that for all $B_{ij}$ the pattern $\sigma(B_{ij})$ has no crossings. Let $j^*$ be the maximal height of the bricks. Then any vertical segment in any  pattern $\sigma^n(B_{ij})$ has length at most $2j^*(\lambda_2-1)$.
\end{proposition}

\emph{Proof:} Let $[(x,y),(x,y')]$ be a vertical segment. The points $(x,y)$ and $(x,y')$  lie in some $\sigma(B_{ij})$ and $\sigma(B_{i'\!j'})$. If these are the same, then $|y-y'|<j^*(\lambda_2-1)$.  If these are different, then the no crossing property implies that they are neighbors of each other, and the length of the segment is at most twice the maximal height minus one of a $\sigma(B_{ij})$, i.e., at most $2j^*(\lambda_2-1)$. \hfill $\square$\\

\noindent
We give an example with $\lambda_1=2$ and $\lambda_2=3$.

\begin{center}
\begin{tikzpicture}[scale=.2,rounded corners=0]
\def\brickoneone(#1){\path (#1) coordinate (P0); \draw [fill=red!40!yellow, draw=black!40, very thick]  (P0) rectangle +(1,1)}
\def\bricktwoone(#1){\path (#1) coordinate (P0);  \draw [fill=red!20!brown, draw=black!40, very thick] (P0) rectangle +(2,1)}
\def\sigma(#1,#2){%
\path (#2) coordinate (P0);
\path (P0)++(-1,0) coordinate (P110);\path (P0)++(0,1) coordinate  (P111);
\path (P0)++(1,1) coordinate  (P112);\path (P0)++(0,2) coordinate  (P113);
\ifthenelse{#1=11}{\bricktwoone(P110);\brickoneone(P111); \brickoneone(P112);\bricktwoone(P113);} {};
\path (P0)++(-1,0) coordinate (P210);\path (P0)++(1,0) coordinate  (P211);
\path (P0)++(2,0) coordinate  (P212);\path (P0)++(0,1) coordinate  (P213);
\path (P0)++(1,1) coordinate  (P214);\path (P0)++(3,1) coordinate  (P215);
\path (P0)++(0,2) coordinate  (P216);\path (P0)++(2,2) coordinate  (P217);
\ifthenelse{#1=21}{\bricktwoone(P210);\brickoneone(P211);\brickoneone(P212);\brickoneone(P213);
\bricktwoone(P214);\brickoneone(P215); \bricktwoone(P216);\bricktwoone(P217);}{};
};
\path (0,0) coordinate (B11); \path (7,0) coordinate (SB11); \path (20,0) coordinate (B21); \path (29,0) coordinate (SB21);
\brickoneone(B11); \sigma(11,SB11); \node at (4,0.4) {$\mapsto$};
\bricktwoone(B21); \sigma(21,SB21); \node at (24.7,0.4) {$\mapsto$};
\end{tikzpicture}
\end{center}

\begin{figure}[h!]%
\begin{center}
 \includegraphics[height=7cm]{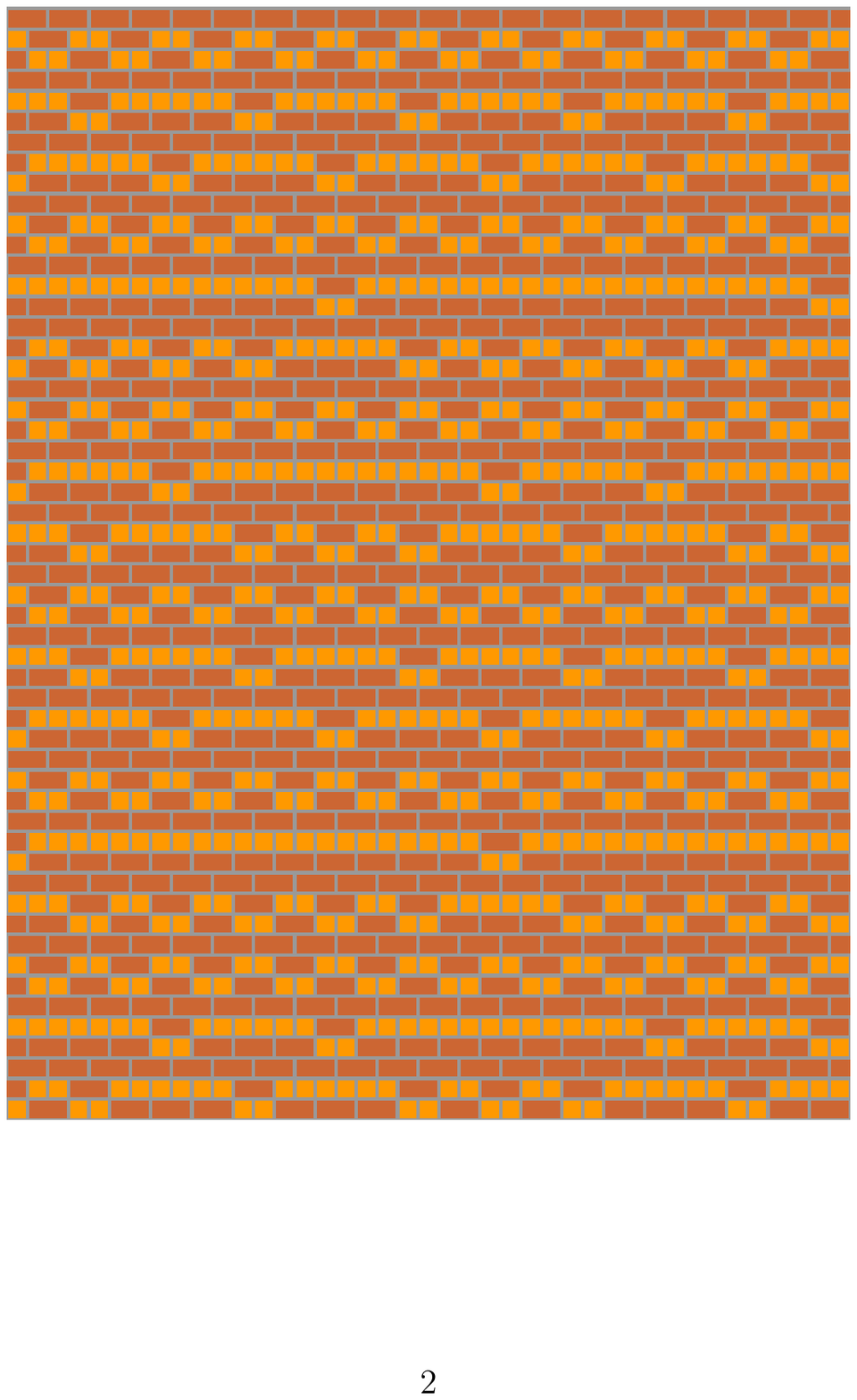}\\  %
  \caption{\footnotesize A clip from $\sigma^5(B_{21})$.}\label{Fig:clip}
\end{center}
\end{figure}

This substitution is of PP-type, and has no crossings. According to Proposition~\ref{Prop:nocross} (with $j^*=1$), the length of vertical segments in any pattern generated by this substitution is at most 4. Actually a refinement of the proof of Proposition~\ref{Prop:nocross} shows that is even bounded by 3, and this bound is achieved as can be verified in Figure~\ref{Fig:clip}.

We end this section with a few remarks on the \emph{number} of bricks that are needed in a pattern. To this end one introduces the matrix $M_\sigma$, with
entries $M_\sigma(ij,i'\!j')$ which by definition are equal to the number of times the brick $B_{i'\!j'}$ occurs in the pattern $\sigma(B_{ij})$, for all $B_{ij}$ and $B_{i'\!j'}\in \mathcal{A}$ (in the literature sometimes the transpose of $M_\sigma$ is taken as definition).

\begin{proposition} The maximal eigenvalue of $M_\sigma$ is $\lambda_1\lambda_2$. The asymptotic frequencies of the bricks are given by the normalized left eigenvector corresponding to this eigenvalue.
\end{proposition}

\emph{Proof:} This proposition is well-known, see e.g.~\cite{Solom}. Here we give a short proof, in a slightly more general situation. We will assume that just as in the self-similar case $\sigma$ satisfies for all $B_{ij}\in\mathcal{A}$
$$ \Ar(\sigma(B_{ij}))= \lambda_1\lambda_2\,\Ar(B_{ij}). $$
Note that our PP-type substitutions satisfy this equation.
Now let $v$ be the vector with entries $\Ar(B_{ij})$, and $w$ the vector with entries $\Ar(\sigma(B_{ij}))$. Then, since $\Ar(\sigma(B_{ij}))$ is the sum of the areas of the $B_{i'\!j'}$ times the number of times $B_{i'\!j'}$ occurs in $\sigma(B_{ij})$, we obtain
$$M_\sigma v= w= \lambda_1\lambda_2 v.$$
Conclusion: $v$ is a right eigenvector of $M_\sigma$ with the eigenvalue $\lambda_1\lambda_2$. But since all the entries of $v$ are strictly positive,  $\lambda_1\lambda_2$
has to be the Perron-Frobenius eigenvalue of $M_\sigma$. The second statement in the proposition is a  well-known result (see e.g.~\cite{Queff}), generalized by Peyri\`ere (see \cite{Pey}) to  higher dimensions:  the normalized left eigenvector of $M_\sigma$ corresponding to the Perron-Frobenius eigenvalue equals the asymptotic frequency of the bricks in large patterns. \hfill $\square$

As a   final remark we mention that the PP-property can be extended in the obvious way to substitutions that in addition have a periodic parallelogram structure in the horizontal direction. We illustrate this with an example:
\begin{center}
 \includegraphics[height=1.5cm]{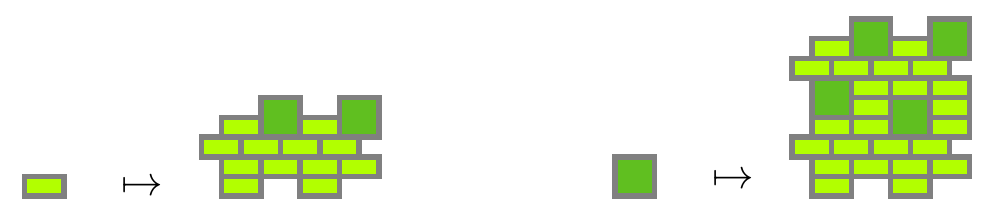}
\end{center}
The expansions are here both 4, there are no crossings, and it is easily shown that $v_{\max}=3$. The patterns $\sigma^2(B_{21})$ and $\sigma^2(B_{22})$
are given by

\begin{center}
\includegraphics[height=4cm]{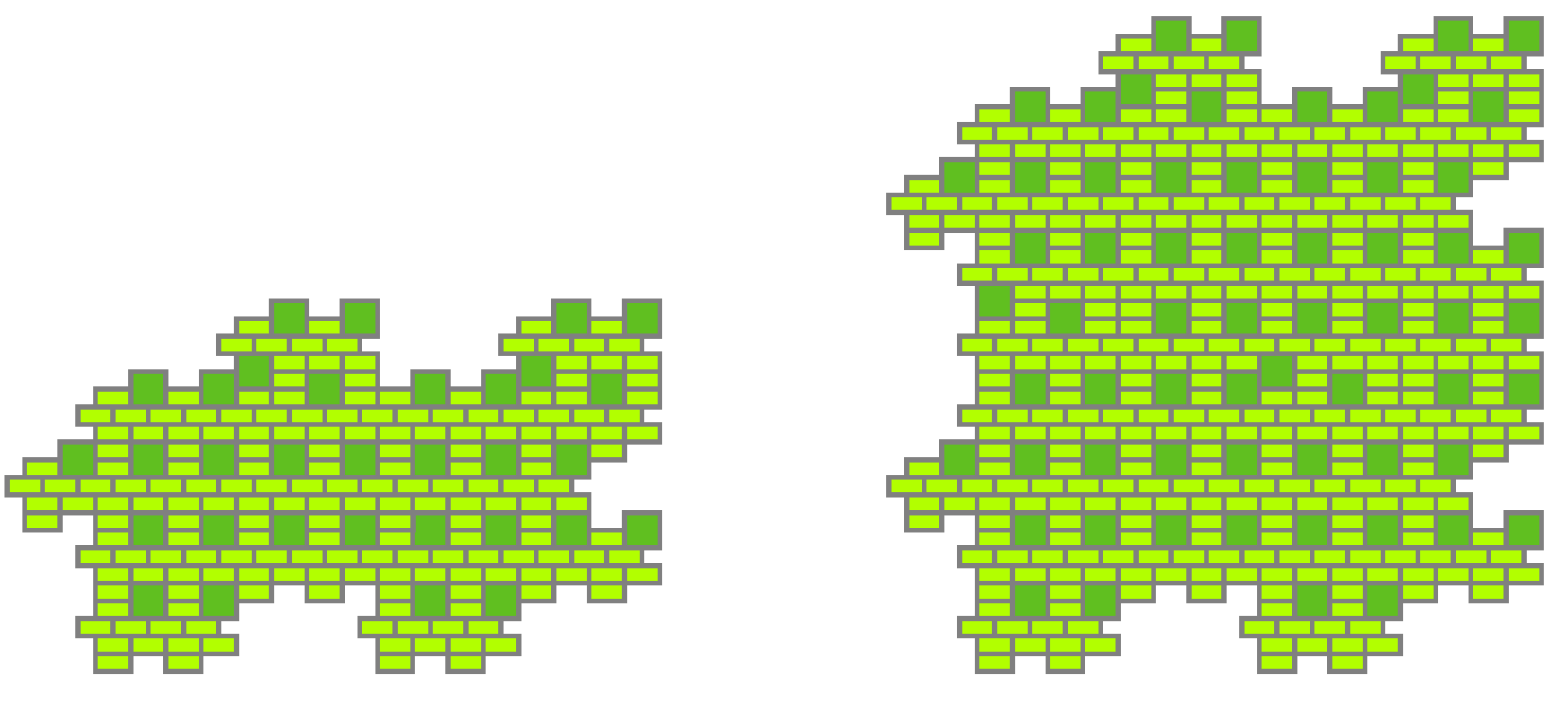}%
\end{center}

\section{Random brick wall patterns}\label{sec:tiles}

 The interpretation of the requirements [R1] (absence of regularity) and [R4] (pleasing aesthetics) is of course for a large part a matter of taste, and some might
consider the pattern in Figure~\ref{Fig:clip} too `rigid'. One way of making the patterns less rigid is to introduce randomness in the construction. As  an introduction we start with a very simple self-similar example, which hence does not satisfy requirement [R3].
We take the $1\times2$ and $2\times2$ bricks with random substitution rules given by

\begin{center}
\begin{tikzpicture}[scale=.2,rounded corners=0]
\def\bricktwoone(#1){\path (#1) coordinate (P0); \draw [fill=red!40!yellow, draw=black!40, very thick]  (P0) rectangle +(1,2)}
\def\bricktwotwo(#1){\path (#1) coordinate (P0);  \draw [fill=red!20!brown, draw=black!40, very thick] (P0) rectangle +(2,2)}
\def\sigma(#1,#2){
\path (#2) coordinate (P0);
\path (P0)++(0,0) coordinate (P1200); \path (P0)++(0,2) coordinate (P1210); \path (P0)++(1,2) coordinate (P1220);
\path (P0)++(0,2) coordinate (P1202); \path (P0)++(0,0) coordinate (P1212); \path (P0)++(1,0) coordinate (P1222);
\ifthenelse{#1=10}{\bricktwotwo(P1200);\bricktwoone(P1210); \bricktwoone(P1220)} {};
\ifthenelse{#1=12}{\bricktwotwo(P1202);\bricktwoone(P1212); \bricktwoone(P1222)} {};
\path (P0)++(0,0) coordinate (P220); \path (P0)++(1,0) coordinate (P221); \path (P0)++(3,0) coordinate (P222);
\path (P0)++(0,2) coordinate (P2230); \path (P0)++(2,2) coordinate (P2240); \path (P0)++(3,2) coordinate (P2250);
\path (P0)++(2,2) coordinate (P2232); \path (P0)++(0,2) coordinate (P2242); \path (P0)++(1,2) coordinate (P2252);
\ifthenelse{#1=20}{\bricktwoone(P220);\bricktwotwo(P221);\bricktwoone(P222);\bricktwotwo(P2230);\bricktwoone(P2240);\bricktwoone(P2250)} {};
\ifthenelse{#1=22}{\bricktwoone(P220);\bricktwotwo(P221);\bricktwoone(P222);\bricktwotwo(P2232);\bricktwoone(P2242);\bricktwoone(P2252)} {};
};
\path (0,0) coordinate (B1); \path (8,3) coordinate (SB10); \path (8,-5) coordinate (SB12);
\path (22,0) coordinate (B2); \path (30,3) coordinate (SB20); \path (30,-5) coordinate (SB22);
\bricktwoone(B1); \draw [->] (3,2.5) -- (6,5); \sigma(10,SB10); \draw [->] (3,-0.5) -- (6,-3); \sigma(12,SB12);
\bricktwotwo(B2); \draw [->] (25,2.5) -- (28,5); \sigma(20,SB20); \draw [->] (25,-0.5) -- (28,-3); \sigma(22,SB22);
\end{tikzpicture}
\end{center}

\begin{figure}[h!]%
\begin{center}
\includegraphics[height=4cm]{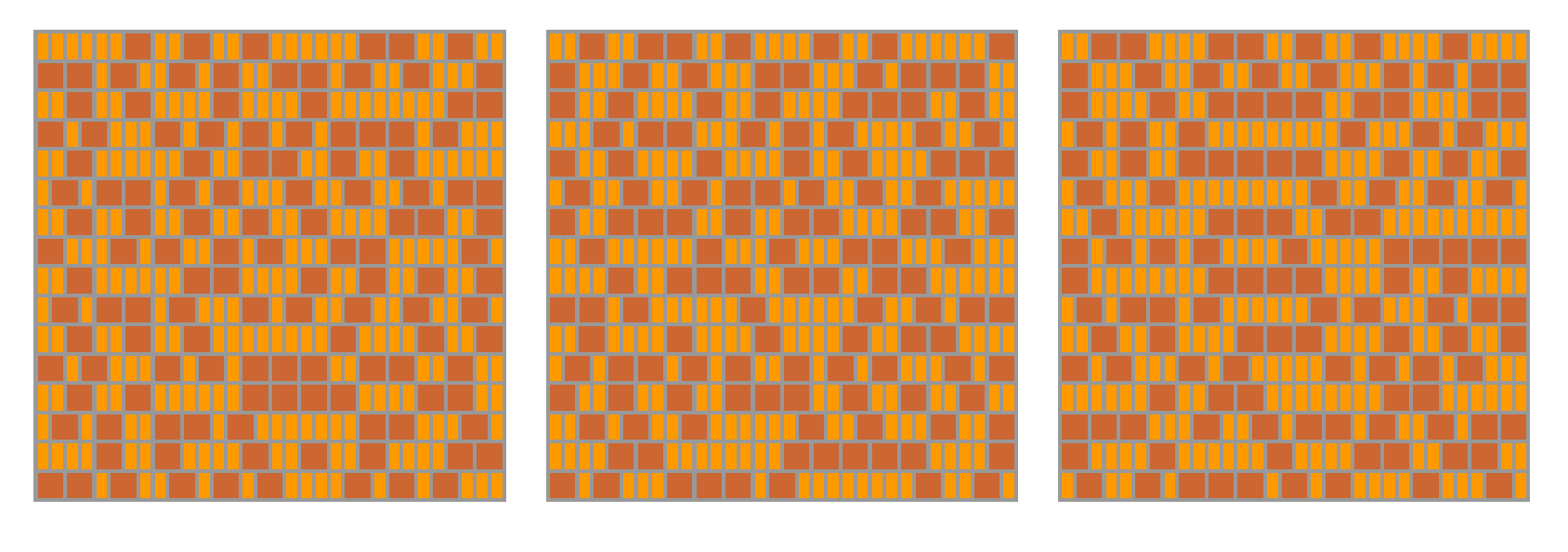}%
  \caption{\footnotesize Three realizations of $\sigma^4(B_{22})$.}\label{Fig:random-4}
\end{center}
\end{figure}

Here random means that each expanded brick is replaced by one of the two possibilities with probability $1/2$, independently of the other replacements.
In this simple example the substitution matrix $ M_{\sigma}$ is non-random, since the two possibilities are just permutations of each other.
A simple computation learns that for $n=1,2,\dots$
$$ M_{\sigma}^n=\left( \begin{array}{cc}
2\cdot 4^{n-1} &  4^{n-1} \\
4^{n} & 2\cdot 4^{n-1}  \end{array} \right).
$$
Note that different sequences of coin flips will lead to different patterns, so there will be
$$2^{\,1+6+ \dots + 6\cdot 4^{n-1}}= 2^{\,2\cdot 4^n-1}$$
different realizations of $\sigma^{n+1}(B_{22})$, since the total number of bricks in $\sigma^{n}(B_{22})$ is equal to the sum
of the two entries in the lower row of $ M_{\sigma}^n$. This means that the three realizations in Figure~\ref{Fig:random-4} are picked from a set of $2^{\,2\cdot 4^3-1}=2^{127}\approx 10^{39}$ elements!

In general the matrix $ M_{\sigma}$ will be random, and its expectation will give us the mean growth of the number of bricks in the
patterns $\sigma^n(B_{ij})$. The theory of multi-type branching processes can be applied to obtain further information on these patterns.

We next consider a pseudo-self-similar random substitution, where the coin is flipped with a success probability $p$ for some $p\in[0,1]$:

\begin{center}
\begin{tikzpicture}[scale=.2,rounded corners=0]
\def\brickonetwo(#1){\path (#1) coordinate (P0); \draw [fill=blue!30!white, draw=black!40, very thick]  (P0) rectangle +(1,2)}
\def\bricktwotwo(#1){\path (#1) coordinate (P0);  \draw [fill=blue!15!gray, draw=black!40, very thick] (P0) rectangle +(2,2)}
\def\sigma(#1,#2){
\path (#2) coordinate (P0);
\path (P0)++(-1,0) coordinate (P1200);         \path (P0)++(-1,0) coordinate (P1201);
\path (P0)++(0,0) coordinate (P1230);
\path (P0)++(0,2) coordinate (P1210);          \path (P0)++(0,2)  coordinate (P1211);
\path (P0)++(1,2) coordinate (P1220);
\ifthenelse{#1=120}{\brickonetwo(P1200);\brickonetwo(P1230);\brickonetwo(P1210); \brickonetwo(P1220)}{};
\ifthenelse{#1=121}{\bricktwotwo(P1201); \bricktwotwo(P1211)} {};
\path (P0)++(-1,0) coordinate (P2200);         \path (P0)++(-1,0) coordinate (P2201);
\path (P0)++(0,0) coordinate (P2270);
\path (P0)++(1,0) coordinate (P2210);          \path (P0)++(1,0)  coordinate (P2211);
\path (P0)++(2,0) coordinate (P2220);          \path (P0)++(0,2)  coordinate (P2221);
\path (P0)++(0,2) coordinate (P2230);          \path (P0)++(2,2)  coordinate (P2231);
\path (P0)++(1,2) coordinate (P2240);
\path (P0)++(2,2) coordinate (P2250);
\path (P0)++(3,2) coordinate (P2260);
\ifthenelse{#1=220}{\brickonetwo(P2200);\brickonetwo(P2270);\brickonetwo(P2210);\brickonetwo(P2220); \brickonetwo(P2230);\brickonetwo(P2240);\brickonetwo(P2250);\brickonetwo(P2260)} {};
\ifthenelse{#1=221}{\bricktwotwo(P2201);\bricktwotwo(P2211);\bricktwotwo(P2221);\bricktwotwo(P2231)} {};
};
\path (0,0) coordinate (B1); \path (10,3) coordinate (SB120); \path (10,-5) coordinate (SB121);
\path (22,0) coordinate (B2); \path (32,3) coordinate (SB220); \path (32,-5) coordinate (SB221);
\brickonetwo(B1); \draw [->] (3,2.5) -- (6,5); \sigma(120,SB120); \draw [->] (3,-0.5) -- (6,-3); \sigma(121,SB121);
\node [above] at (3.2,3.75) {\footnotesize $1-p$}; \node [above] at (4.5,-1.75) {\footnotesize $p$};
\bricktwotwo(B2); \draw [->] (25,2.5) -- (28,5); \sigma(220,SB220); \draw [->] (25,-0.5) -- (28,-3); \sigma(221,SB221);
\node [above] at (25.2,3.75) {\footnotesize $1-p$}; \node [above] at (26.5,-1.75) {\footnotesize $p$};
\end{tikzpicture}
\end{center}

The maximal length of a vertical segment now is a random variable $V_{\rm max}$, and in this pseudo-self-similar example the behavior
 of $V_{\rm max}(p)$  as a function of $p$ is particularly interesting. Note that $V_{\rm max}(0)\equiv\infty$, but that $V_{\rm max}(1)\equiv 2$.

\begin{figure}[h!]%
\begin{center}
\includegraphics[height=2.8cm]{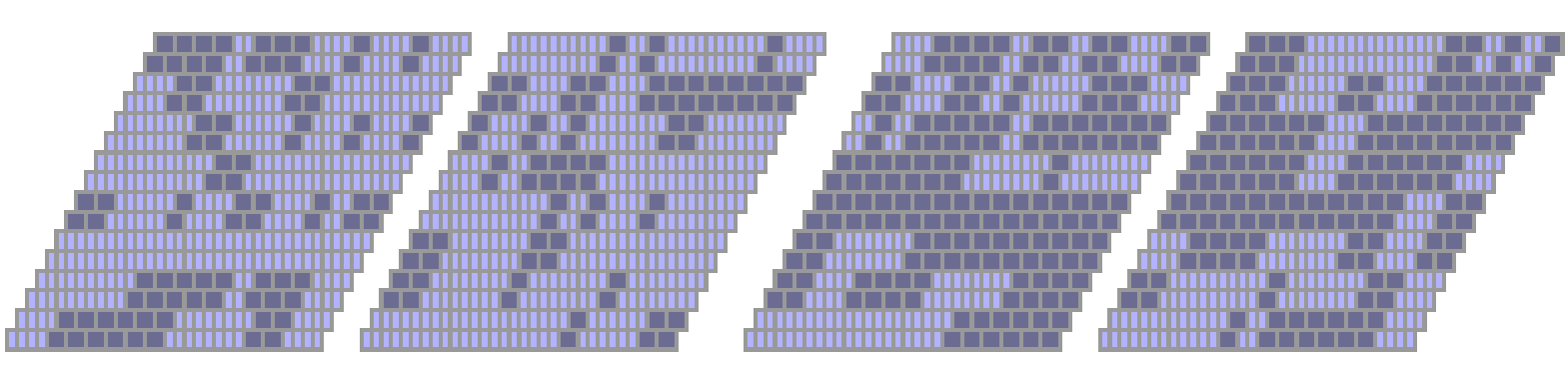}%
  \caption{\footnotesize Left: two realizations of $\sigma^4(B_{22})$ with $p=1/3$. Right: idem with $p=2/3$.}\label{perco}
\end{center}
\end{figure}

 Figure~\ref{perco} shows 2 realizations for $p=1/3$ and $p=2/3$ each.

\end{document}